\documentclass[11pt]{amsart}
\usepackage{amsfonts,amssymb,amsmath,eucal}
\usepackage[all]{xy}
\begin{document}

\newcommand{\R}{{\mathbb R}}
\newcommand{\G}{{\mathbb G}}
\newcommand{\U}{{\rm U}}
\newcommand{\SO}{{\rm SO}}
\newcommand{\mH}{{\mathbb H}}
\newcommand{\SU}{{\rm SU}}
\newcommand{\Oh}{{\rm O}}
\newcommand{\T}{{\mathbb T}} 
\newcommand{\hT}{{\widehat{T}}}
\newcommand{\Z}{{\mathbb Z}}
\newcommand{\Q}{{\mathbb Q}}
\newcommand{\hG}{{\widehat{G}}}
\newcommand{\tK}{{\widetilde{K}}}
\newcommand{\C}{{\mathbb C}}
\newcommand{\Sp}{{\rm Sp}}
\newcommand{\Sl}{{\rm Sl}}
\newcommand{\F}{{\mathbb F}}
\newcommand{\fr}{{\rm free}} 
\newcommand{\op}{{\rm op}}
\newcommand{\pt}{{\rm pt}}
\newcommand{\res}{{\rm res}}
\newcommand{\tc}{{\rm TC}}
\newcommand{\geo}{{\rm geo}}
\newcommand{\hot}{{\rm hot}}
\newcommand{\li}{{\rm li}}
\newcommand{\Spec}{{\rm Spec}}

\title{Geometric Tate-Swan cohomology of equivariant spectra}
\author{Jack Morava}
\address{The Johns Hopkins University,
Baltimore, Maryland 21218}
\email{jack@math.jhu.edu}
\date{20 November 2012}
\begin{abstract}{We sketch a quick and dirty geometric approach to
the Tate-Swan cohomology of equivariant spectra, illustrating it
with conjectural applications to Atiyah-Segal $K$-theory of circle
actions, and a possible geometric model for the topological cyclic
homology of the sphere spectrum.}\end{abstract}

\maketitle \bigskip

\noindent{\it For Michael and Graeme}\bigskip

\section{Introduction}\bigskip

\noindent
These variations on the theme of Tate-Swsn cohomology of spectra 
are very old-fashioned, but my hope is to illustrate the flexibility and
naturality of geometric methods. The first section below summarizes the 
general framework, while the second applies those ideas to Atiyah-Segal 
equivariant cohomology of circle actions. The third section proposes a model 
for the topological cyclic homology of the sphere spectrum.\bigskip

\noindent
It is a pleasure to acknowledge helpful conversations with Bjorn Dundas and Mike Hill
about the material in the first section, and with Lars Hesselholt about the third.
The howlers, however, are my own inventions. \bigskip

\noindent
{\bf 1.1} Let $G$ be a (locally connected) $d$-dimensional compact Lie group, and let $E_G$ be
multiplicative cohomology theory, which admits a reasonable interpretation as a $G$-equivariant
cobordism theory of manifolds with `$E$-structure' of some sort: the underlying non-equivariant
theory could be, for example, the sphere spectrum, the noncommutative spectrum $M\xi$ of [2],
$M\U, \; M\Oh$, or (depending on how hard we're willing to work [1]) $H\Z$. It seems likely that 
some version of the work of Baum, Douglas [3] and others can also be handled by these methods. 
\bigskip

\noindent
We consider compact $E$-oriented manifolds $Z$ with compatible $G$-structure, with boundary
\[
\partial Z \; = \; \partial_0 Z \cup \partial_\fr Z
\]
partitioned into transversally intersecting parts, such that the $G$-action on $\partial_\fr Z$ is in
fact free. I'll call such a manifold {\bf closed} if $\partial_0 Z = \emptyset$: in other words if the
$G$-action on $Z$ is free on the entire boundary. Two such manifolds will be said to be
cobordant if there is a manifold $W$ in this class, together with a transversal partition 
\[
\partial W \; = \; (Z_+ \sqcup Z^\op_-) \cup V
\]
with $G$-action free on $V$, and
\[
V \cap (Z_+ \sqcup Z^\op_-) \; = \; \partial Z_+ \sqcup \partial Z^\op_- \;;
\]
this is all completely classical [6 I \S 4]. I'll write (some variant of) $t^G_nE$ for the abelian
group of cobordism classes of such manifolds with group action free on the boundary, or more generally 
$t^G_*E(X)$ for the graded group of equivalence classes of manifolds mapped equivariantly to a pointed
$G$-space $X$; this defines a $G$-equivariant geometric homology theory. Stabilizing (ie taking 
the direct limit over suspensions from a suitable universe of $G$-representations [5]; in alternate 
terms, sheafifying a presheaf) defines, if we're lucky, a theory represented by a suitable 
$G$-spectrum\begin{footnote}{In particular, if we're not plagued by phantoms. I'm omitting details
because the general constructions of [12] are quite convenient for these technical issues.}\end{footnote}. 
\bigskip

\noindent
{\bf 1.2} For reasonable $E$-structures, the Cartesian product $Z \times Z'$ of two closed $G$-manifolds,
given the diagonal $G$-action, can be smoothed to be another; more precisely, the boundary
\[
\partial (Z \times Z') \; = \; \partial Z \times Z' \cup_{\partial Z \times \partial Z'} Z \times \partial
Z'
\]
has a natural smooth structure [6 I \S 3] compatible with the $G$-action. This defines an
external product
\[
t^G_nE(X) \otimes t^G_mE(Y) \to t^G_{n+m}E(X \wedge Y)
\]
and, in particular, makes $t^G_*E(S^0)$ into a graded ring. \bigskip

\noindent
{\bf 1.3 Examples}\bigskip

\noindent
1) The interval $[-1,+1]$, regarded as an unoriented manifold with $\Z_2$-action defined by $x
\mapsto -x$, defines a class [4, 16-18]
\[
w^{-1} \in t^{\Z_2}_1 M\Oh(S^0) \;.
\]
2) The closed unit disk $D \subset \C$, regarded as a complex-oriented manifold with action
\[
u,z \mapsto uz : \T \times D \to D
\]
of the unit circle $\T = \partial D$, defines a class $c^{-1} \in t^\T_2 M\U(S^0)$.\bigskip

\noindent
3) The closed unit ball $B(\mH)$ in the quaternions, regarded as an $\mH$-oriented manifold
with action 
\[
\SU(2) \times B(\mH) \to B(\mH)
\]
by multiplication of the unit quaternions, defines a class $\wp^{-1} \in t^{\SU(2)}_4M\Sp(S^0)$. \bigskip

\noindent
{\bf 1.4} Constructions of Quillen (again, in good situations)  associate to an $n$-manifold $Z$ with
{\bf empty} boundary, an $E$-structure and a $G$-action, the class
\[
[Z \times_G EG \to S^0 \wedge_G EG_+] \in E^{-n}_G(S^0)
\]
of its homotopy quotient. Similarly, if the action of $G$ on $Z$ is in fact free, its geometric
quotient
\[
[Z/G \to BG_+] \in E_{n-d}(BG_+)
\] 
identifies $Z$ as a principal $G$-bundle. Together these constructions fit into the fundamental exact 
sequence 
\[ 
\xymatrix{
\dots \ar[r]^\phi & E^G_n(S^0) \ar[r]^\rho & t^G_nE(S^0) \ar[r]^<<<<\partial & E_{n-d-1}(BG_+) \ar[r] & \dots 
}
\]
of $E_G^*$-modules, with $\rho$ a ring homomorphisms, and $\phi$ the forgetful map from free to 
unrestricted group actions. [The product of a general $G$-manifold and a manifold with free $G$-action,
given the diagonal group action, is a manifold with free action.]\bigskip

\noindent
{\bf 1.5.1 Examples}, continued:
\[
\begin{array}{ll}
t^*_{\Z_2}H\Z_2(S^0)     &    =  \; \;  \Z_2[w^{\pm 1}] \\
t^*_\T H\Z(S^0)          &    =  \; \;  \Z[c^{\pm 1}] \\
t^*_{\SU(2)} H\Z(S^0)    &    =  \; \;  \Z[\wp^{\pm 1}]
\end{array}
\]
(where $w^{-1}, \; c^{-1}$ and $\wp^{-1}$ are the images of the corresponding cobordism classes
under the Steenrod cycle map).\bigskip

\noindent
In particular, in 
\[
\dots \to H^{-*}(B\T_+,\Z) = \Z[c] \to t^{-*}_\T H\Z(S^0) = \Z[c^{\pm 1}] \to H_{*-2}(B\T_+,\Z) =
\Z[\gamma_n \:|\: n \geq 0] \to \dots \;,
\]
the product $c^{-n}$ represents the class of the unit ball in $\C^n$, with $\T$ acting as
multiplication. The boundary map sends it to the divided power
\[
\partial c^{-n} = [S^{2n-1}/\T = \C P^{n-1} \subset \C P^\infty = B\T] = \gamma_{n-1} \;,
\]
Kronecker dual to the $(n-1)$st power of the usual first Chern class $c \in H^2(B\T,\Z)$. The
classes $w^{-1}$ and $\wp^{-1}$ are similarly related to the usual first Stiefel-Whitney and
Pontrjagin classes.\bigskip

\noindent
{\bf 1.5.2} More generally, $t^*_\T M\U(S^0)$ is a formal Laurent series ring $M\U^*((c))$, and 
$M\U_*B\T$ is a free $M\U_*$ - module on generators $\beta_n$ Kronecker
dual to the Chern classes $c^n \in M\U^{2n}(B\T)$, satisfying 
\[
\beta(s_0 +_F s_1) = \beta(s_0)\beta(s_1) \;,
\]
where $\beta(s) = \sum \beta_n s^n$ and $s_0 +_F s_1 = F_{M\U}(s_0,s_1)$ is the universal formal
group law [19]; the argument above generalizes, implying $\partial c^{-n} = \beta_{n-1}$.\bigskip

\noindent
The map $\epsilon : B\T \to \pt$ defines a kind of residue homomorphism
\[
\xymatrix{
\res : t^*_\T M\U(S^0) \ar[r]^\partial & M\U_{-*-2}B\T_+ \ar[r]^{\epsilon_*} & M\U_{-*-2}(S^0)
}
\]
satisfying $\res \; (c^{-n}) = \delta_{n,1}$; this lets us write
\[
(c^n,\beta_m) = \res \; (c^n \cdot c^{-m-1}) \; \; (n,m \geq 0)
\]
for the Kronecker product, ie 
\[
\beta_m = \res \; (c^{-m-1} \cdot -) \in M\U_{2m}B\T_+ \to {\rm Hom}^{-2m}_{M\U}(M\U^*(B\T_+),M\U(S^0)) \;.
\] 
Thus if $f \in t^*_\T M\U(S^0), g \in M\U^*(B\T_+)$, we have [17]
\[
(\partial f)(g) = \res \; (f \cdot g) \in M\U_*(S^0) \;.
\]
{\bf 1.5.3} The Segal conjecture for finite groups [15] supplies another class of examples: 
after a suitable completion, the exact sequence above simplifies to an equivalence
\[
t_{\hat{G}} \; \sim \; \vee_{e \neq H < G} \widehat{BW}_H 
\]
where $W_H = N(H)/H$ is a kind of Weyl group. \bigskip

\section{Classical $K_\T$}\bigskip

\noindent
Rudimentary knowledge of equivariant homotopy theory tells us that an ordinary
(nonequivariant) cohomology theory can have more than one equivariant extension, and the
account above ignores this. Rather than confront that issue, I'll consider two examples related
to this question, which I still do not understand well enough. \bigskip

\noindent
{\bf 2.1} If $G$ is a finite group, G. Wilson's identification [25 Prop 1.2] of $\tilde{K}_*(BG)$ 
as the $\Q/\Z$-dual of the augmentation ideal of the completed representation ring $\hat{R}(G)$, concentrated in
odd degree, identifies [11, 12] the Tate cohomology
\[
t_GK^\hot(S^0) \; = \; \Z \oplus \hat{R}_+(G) \otimes \Q
\]
of homotopy-theoretic $G$-equivariant $K$-theory (cf also [3]). It seems plausible that something
similar holds for Atiyah-Segal equivariant $K$-theory, but I do not know of a proof. 
In any case Tate cohomology defines an interesting analog
\[
K^\hot_G \to t_G K^\hot
\]
of the Chern character. \bigskip

\noindent
{\bf 2.2} The representation ring $R(\T)$ of the circle is the Atiyah-Segal $\T$-equivariant
$K$-theory $K[\chi^{\pm 1}]$ of a point ($K = \Z$ as a $\Z_2$-graded ring, and $\chi = \exp(2\pi 
i\theta)$ is the standard character $\T = \R/\Z \to \C$), while 
\[
K^\hot_\T(S^0) := K(B\T_+) \cong \Z[[t]]
\]
with $\chi = 1 - t$. The homology
\[
K_*(B\T_+) \; = \; K[b_n \:|\; n\geq 0]
\]
is the algebra (under Pontrjagin product) generated by elements satisfying the identity
\[
b(s_0) b(s_1) = b(s_0 + s_1 - s_0s_1) \;,
\]
where 
\[
b(s) = (1-s)^t = \sum_{n \geq 0} b_n s^n = \sum {t \choose n} (-s)^n 
\]
so the exact sequence 
\[
0 \to K(B\T_+) = K[[t]] \to t_\T K^\hot(S^0) = \Z((t)) \to \Z[b_n \:|\; n \geq 0] \to 0
\]
is much like Ex. 1.5.2. \bigskip

\noindent
{\bf 2.3} The homomorphism $\chi \mapsto 1 - t$ defines a completion
\[
t_\T K(S^0) \cong R(\T)[(1-\chi)^{-1}] \cong K[\chi^{\pm 1},(1-\chi)^{-1}] \to t_\T K^\hot(S^0) = K((t))
\]
at the identity $\chi = 1$ of the multiplicative group Spec $K_\T = \G_m$ (cf [22]), analogous to a point
\[
\Spec \; K((t)) \to \Spec \; t_\T K(S^0) \; \sim \; {\mathbb P}^1 - \{0,1,\infty\} \;.
\]
The group
\[
\Sigma_3 \cong \{\sigma,\tau \:|\: \tau^2 = 1, \; \sigma^3 = 1, \; \tau^{-1} \sigma \tau = \sigma^{-1} \}
\]
acts by fractional linear transformations $\tau(\chi) = \chi^{-1}, \; \sigma(\chi) = (1-\chi)^{-1}$ on 
the projective line over $K$, permuting the points $\{0,1,\infty\}$. I don't know if this lifts to any 
kind of action on the functor $t_\T K$.\bigskip

\noindent
Writing $c = 2\pi i \theta$ extends the usual Chern character
\[
K(B\T_+) \to \Q[[c]]
\]
to a specialization $\Spec \; \Q((c)) \to \Spec \; t_\T K(S^0)$ of the map above, which sends $q := (1 - \chi^{-1})^{-1}$
to $c^{-1}(1 + \dots) \in c^{-1}\Q[[c]]$; while $q \mapsto \exp(\hbar)$ similarly maps $\chi$ to $\hbar^{-1}\Q[[\hbar]]$,
defining another formal point 
\[
\Spec \; \Q((\hbar)) \to  \Spec \; t_\T K(S^0) \;.
\]
The first of these corresponds to completion near $c = 0$ (ie $\chi = 1$), while the second is completion 
near $\chi = \infty$. In terms of the distribution 
\[
\li_1(x) = \log|1 - e^x| \; (\; \equiv - \log|x| \; {\rm mod \; smooth \; functions \; of} \; x)
\]
defined by the composition
\[
\log \circ \left[\begin{array}{cc}
                -1 & 1 \\
                 0 & 1 
                 \end{array}\right] \circ \exp
\]
[9 \S 3 Cor 2] we have
\[
- \hbar = \li_1(-c), \; - c = \li_1(-\hbar) \;,
\]
suggesting that $\hbar \to 0 \Longleftrightarrow c \to \infty$ is a kind of semiclassical limit \dots \bigskip 
  
\noindent
Note finally that $\Z[\chi^{\pm 1},(1-\chi)^{-1}]$ has Tate's Laurent series ring $\Z((\chi))$ (which accomodates $K$-theory 
with the $\sigma$-orientation) as a completion. This suggests the interest of the corresponding completion of $K_\T(LX)$ (or of
$t_\T K(LX)$) as a model for elliptic cohomology [13]. On the other hand, a theorem of Goodwillie [10] suggests that 
$t_\T K^\hot(LK)$ (and hence $K_\T (LX)\otimes_{K_\T} K\Q((c))$) sees only the fundamental group of $X$ : the left vertical arrow
in the diagram
\[
\xymatrix{
K^\hot_\T(-)[c^{-1}] \ar[d] \ar[dr] \ar[r]^{ch} & H^*_\T(-,K\Q)[c^{-1}] \\
K_\T(-)((c)) \ar[r] & K^\hot_\T(-)[c^{-1}]\otimes \Q \ar@{.>}[u] }
\]
is an equivalence of functors, while the right vertical arrow is induced
by the injective transformation
\[
M \to M \otimes_{\Z((c)) \otimes \Q} \Q((c)) \;.
\]
\bigskip \bigskip

\section{A model for $\tc(S^0)$} \bigskip

\noindent
Work of B\"okstedt, Hsiang, and Madsen identifies the topological cyclic homology
$\tc(S^0)$ of the sphere, after $p$-adic completion, with a similar completion of
the spectrum $S^0 \vee \Sigma \C P^\infty_{-1}$. This section records the construction 
of a geometric model $\tc^\geo$ for the latter object, which has an interesting multiplicative
structure. However we make no attempt to construct a map to or from $\tc(S^0)$ itself.\bigskip 

\noindent
{\bf 3.1} A reasonable cobordism theory $\Omega^G_*$ has an associated cobordism theory 
$\Omega^G_* \oplus \Omega^G_{*-1}$ of $G$-manifolds with boundary [6 I \S 4], 
represented by $(S^0 \vee S^1) \wedge MG$; where $MG$ is the spectrum representing 
$\Omega^G_*$. The homotopy fiber $\tc^\geo$ of the composition
\[
\Sigma B\T_+ \to S^0 \to S^0 \vee S^1 = \T_+
\]
(defined by the stable circle-transfer [14], followed by the obvious inclusion) can be 
interpreted as representing a cobordism theory of framed manifolds $M$ with boundary, together 
with extra data defined by a complex line bundle on the manifold, and a trivialization
of that bundle away from a collar neighborhood of a codimension zero submanifold of 
its boundary: a variation on the Baas-Sullivan theory [1] of cobordism with singularities,
based on framed manifolds in which the boundary $\partial M = \partial_0 M \cup \partial_1 M$ 
is partitioned into two (transversally intersecting) parts, with $\partial_1 M$ carrying
a line bundle trivialized away from $\partial (\partial_1 M)$. The operation $M \mapsto 
\partial_1 M$ thus satisfies $\partial_1 \circ \partial_1 = \emptyset$.\bigskip

\noindent
If $M$ is closed in this relative sense (ie $\partial_1 M = \emptyset$), and 
\[
u,z \mapsto uz : \T \times D \to D
\]
is the usual circle action on the disk $\{z \in \C \:|\: |z = 1 \}$, then (using the
natural framing of the circle bundle $C(\partial_0 M)$)
\[
D \times_\T C(\partial_0 M) \cup_{0 \times \partial_0 M} M \; := \; M_D
\]
is a closed framed manifold. The closed objects in this relative cobordism category are
thus something like the algebraic geometers' varieties bounded by divisors. \bigskip

\noindent
{\bf 3.2} The homotopy exact sequence 
\[
\cdots \to \pi^S_{n+1}(S^0) \oplus \pi^S_n(S^0) \to \tc^\geo_n(S^0) \to \pi^S_{n-1}B\T_+ \to 
\pi^S_n(S^0) \oplus \pi^S_{n-1}(S^0) \to \cdots 
\]
associated to this construction starts at the left by sending a framed $(n+1)$ - manifold with 
boundary to its boundary, regarded as a framed manifold decorated with a trivial complex line bundle. 
The second arrow in the sequence sends a framed $n$ - manifold with boundary and a suitable
complex line bundle, to the class of its boundary, regarded as an element of the 
$(n-1)$ - dimensional bordism group of the classifying space for circle bundles. The circle 
transfer defines the third homomorphism of the sequence, which sends a closed framed 
$(n-1)$ - manifold with a circle bundle over it, to the total space of that bundle (given 
its natural framing), regarded as an $n$ - manifold with empty boundary. \bigskip

\noindent
{\bf 3.3} This cobordism theory has a natural multiplication, defined by the tensor product of
line bundles over the cartesian product of underlying manifolds. $\tc^\geo_0$ is 
generated by a point, and $\tc^\geo_1$ has a rank one part generated by a closed 
interval. $\tc^\geo_{-1}$ has a somewhat unconvincing interpretation as generated by the 
$-1$ - dimensional manifold bounded by the $-2$ - dimensional manifold carrying 
the complex line bundle whose total space is a point \dots \bigskip

\noindent
{\bf 3.4} It follows from the cofibration sequence 
\[
\dots \to K^*(S^0 \vee S^1) \to K^*(\Sigma B\T_+) \to K^*(\tc^\geo) \to \dots
\]
that $K^0(\tc^\geo) \cong \Z$ and
\[
K^1(\tc^\geo) \cong \Z \langle t^k \:|\: k \geq -1 \rangle
\]
with $t = \eta - 1 : B\T \to B{\rm U}$ classifying the Hopf line bundle. In the 1970's 
Segal suggested the homotopy fiber $\Sigma^{-1}K\C^\times$ of the map 
\[
\xymatrix{
K \ar[r]^{2 \pi i} & K\C }
\]
as an interesting model for the algebraic $K$-theory of $\C$ (as a discrete field).
This suggests regarding the image of $t$ in $[\tc^\geo,\Sigma^{-1}K\C^\times]$ as a
topological analog of Borel's regulator [8].

\newpage   

\bibliographystyle{amsplain}

\begin{thebibliography}{99}


\bibitem[1]{1} NA Baas, On bordism theory of manifolds with singularities,
Math. Scand. 33 (1973) 279 - 302

\bibitem[2]{2} A Baker, B Richter, Quasisymmetric functions from a topological
point of view, Math. Scand. 103 (2008) 208 - 242

\bibitem[3]{3} P Baum, RG Douglas. Index theory, bordism, and $K$-homology, in
{\bf Operator algebras and $K$-theory} p. 131, Contemp. Math. 10, AMS
(1982)

\bibitem[4]{4} JM Boardman, Cobordism of involutions revisited, in {\bf Proceedings
of the Second Conference on Compact Transformation Groups} pp. 131 - 151.
Lecture Notes in Math., Vol. 298, Springer (1972)

\bibitem[5]{5} Th. Br\"ocker, EC Hook, Stable equivariant bordism, Math. Z. 129
(1972) 269 - 277

\bibitem[6]{6} PE Conner, EE Floyd, {\bf Differentiable periodic maps}. Ergebnisse der
Mathematik und ihrer Grenzgebiete 33 (1964)

\bibitem[7]{7} T tom Dieck, {\bf Transformation groups}, de Gruyter Studies in 
Mathematics 8 (1987)

\bibitem[8]{8} J Dupont, R Hain, S Zucker, Regulators and characteristic classes of 
flat bundles, in {\bf The arithmetic and geometry of algebraic cycles} 47 - 92,
CRM Proc. Lecture Notes 24, AMS (2000)                              

\bibitem[9]{9}  C Epstein, J Morava, Tempering the polylogarithm, {\tt arXiv:math/0611240} 

\bibitem[10]{10} T Goodwillie, Cyclic homology, derivations, and the free
loopspace, Topology 24 (1985) 187 - 215

\bibitem[11]{11} JPC Greenlees, Tate cohomology in commutative algebra, J. Pure Appl. Algebra
94 (1994) 59 - 83

\bibitem[12]{12} ------, JP May, {\bf Generalized Tate cohomology}, Mem. Amer.
Math. Soc. 113 (1995) no. 543

\bibitem[13]{13} N Kitchloo, J Morava, Thom prospectra for loop group representations,
{\tt arXiv:math/0404541}

\bibitem[14]{14} I Madsen, C Schlichtkrull, The circle transfer and $K$-theory, in
{\bf Geometry and topology: Aarhus (1998)} 307 - 328, Contemp. Math
258, AMS (2000)

\bibitem[15]{15} HR Miller, appendix to DC Ravenel, The Segal conjecture for cyclic
groups and its consequences, Amer. J. Math. 106 (1984) 415 - 446

\bibitem[16]{16} J Morava, Cobordism of involutions revisited, revisited, in {\bf Homotopy
invariant algebraic structures}, Contemp. Math. 239, AMS (1999)

\bibitem[17]{17} ------, Tate cohomology of circle actions as a Heisenberg group,
{\tt arXiv:math/0102132}

\bibitem[18]{18}. ------, Heisenberg groups and algebraic topology, in {\bf Topology, geometry
and quantum field theory} pp. 235 - 246. LMS Lecture Notes 308, Cambridge (2004)

\bibitem[19]{19} D Quillen, Elementary proofs of some results of cobordism theory
using Steenrod operations, Advances in Math. 7 (1971) 29 - 56

\bibitem[20]{20} DC Ravenel, WS Wilson, The Hopf ring for complex cobordism. J. Pure Appl.
Algebra 9 (1976/77) 241 - 280

\bibitem[21]{21} J Rognes, The smooth Whitehead spectrum of a point at odd regular
primes, {\tt arXiv:math/0304384}

\bibitem[22]{22} I Rosu, Equivariant $K$-theory and equivariant cohomology,
{\tt arXiv:math/9912088} 

\bibitem[23]{23} NP Strickland, $K(n)$-local duality for finite groups and
groupoids. Topology 39 (2000) 733 - 772

\bibitem[24]{24} B Walker, Orientations and $p$-adic analysis, available at
{\tt arXiv:0905/0022}

\bibitem[25]{25} G Wilson, $K$-theory invariants for unitary $G$-bordism, Quart. J.
Math. 24 (1973) 499 - 526

\end{thebibliography}

\end{document}